\documentclass[leqno,twoside,a4paper,12pt]{article}
\usepackage{amsmath,amssymb,amsthm}
\usepackage{url}
\usepackage[midshaft,scriptlabels]{diagrams}
\DeclareFontFamily{U}{UWCyr}{}
\DeclareFontShape{U}{UWCyr}{m}{n}{%
  <5> <6> <7> <8> <9>
  <10> <10.95> <12> <14.4> <17.28> <20.74> <24.88> wncyr10
  }{}
\DeclareMathAlphabet{\cyrm}{U}{UWCyr}{m}{n}
\DeclareSymbolFont{cyrm}{U}{UWCyr}{m}{n}
\DeclareSymbolFontAlphabet{\cyrm}{cyrm}
\DeclareMathSymbol{\Evo}{\cyrm}{cyrm}{"03}
\newtheorem{theorem}{Theorem}
\newtheorem{lemma}{Lemma}

\newtheorem{proposition}{Proposition}

\theoremstyle{definition}
\newtheorem{definition}{Definition}

\theoremstyle{remark}
\newtheorem{remark}{Remark}
\DeclareMathOperator{\byd}{\raisebox{-.2ex}{$\overset{\text{\tiny def}}{=}$}}
\DeclareMathOperator{\im}{im}
\DeclareMathOperator{\pro}{pr}
\DeclareMathOperator{\id}{id}

\DeclareMathOperator{\Vol}{Vol}
\DeclareMathOperator{\Imm}{Imm}
\DeclareMathOperator{\Gr}{Gr}
\renewcommand{\to}{\rTo[]\relax}
\newarrow{MapsTo}|--->
\newcommand{\mto}{\rMapsTo[]\relax}
\newarrow{hTo}{boldhook}--->
\newcommand{\hto}{\rhTo[]\relax}

\newcommand{\cprime}{\/{\mathsurround=0pt$'$}}

\newcommand{\bnu}{\boldsymbol{\nu}}
\newcommand{\brh}{\boldsymbol{\rho}}
\newcommand{\bsi}{\boldsymbol{\sigma}}
\newcommand{\bta}{\boldsymbol{\tau}}

\newcommand{\hd}{\bar{d}}
\newcommand{\R}{\mathbb{R}}
\newcommand{\cC}{\mathcal{C}}

\newcommand{\cF}{\mathcal{F}}
\newcommand{\cH}{\mathcal{H}}
\newcommand{\cL}{\mathcal{L}}

\newcommand{\Diff}{\mathrm{Dif{}f}}
\DeclareMathOperator{\alt}{alt}

\DeclareMathOperator{\Hom}{Hom}
\newcommand{\hL}{\Bar{\Lambda}}
\newcommand{\CDiffalt}[2]{\cC\Diff^{\alt}_{(#1)\,#2}}
\newcommand{\CDiff}{\mathcal{C}\mathrm{Diff}}

\newcommand{\odx}{\overline{dx}{}}
\newcommand{\abs}[1]{\lvert#1\rvert}

\newcommand*{\pd}[2]{\mathchoice{\frac{\partial #1}{\partial #2}}
  {\partial #1/\partial #2}{\partial #1/\partial #2} {\partial%
  #1/\partial #2}}
\newcommand{\eval}[2][\right]{\relax
  \ifx#1\right\relax \left.\fi#2#1\rvert}

\let\abs=\envert

\providecommand{\href}[2]{#2}

\newcommand*{\email}[1]{\href{mailto:#1}{\begingroup \urlstyle{rm}\Url{#1}}}
\newcommand*{\eprint}[2][]{\href{http://arXiv.org/abs/#2}%
{\begingroup \Url{arXiv:#2}}}
\title{The geometry of finite order jets of submanifolds\\
and the variational formalism}
\author{Gianni Manno
\\ {\footnotesize Dept.\ of Mathematics, King's College, London}
\\ {\footnotesize email: \email{Gianni.Manno@kcl.ac.uk}}
\and Raffaele Vitolo
\\ {\footnotesize Dept.\ of Mathematics ``E. De Giorgi'', University of Lecce}
\\ {\footnotesize via per Arnesano, 73100 Lecce, Italy}
\\ {\footnotesize and Diffiety Institute, Russia}
\\ {\footnotesize email: \email{Raffaele.Vitolo@unile.it}}
}
\date{}
\begin{document}

\maketitle

\begin{abstract}
  We study the geometry of jets of submanifolds with special interest in the
  relationship with the calculus of variations. We give a new proof of the fact
  that higher order jets of submanifolds are affine bundles; as a by-product we
  obtain a new expression for the associated vector bundles. We use
  Green--Vinogradov formula to provide coordinate expressions for all
  variational forms, \emph{i.e.}, objects in the finite-order variational
  sequence on jets of submanifolds.  Finally, we formulate the variational
  problem in the framework of jets of submanifolds by an intrinsic geometric
  language, and connect it with the variational sequence. Detailed comparison
  with literature is provided throughout the paper.
\end{abstract}

\tableofcontents

\section*{Introduction}

The study of jet spaces was born as the study of the notion of \emph{contact}.
From the analytical viewpoint, $r$-th contact between smooth maps was realized
as the equality of their $r$-th differentials. From the geometrical viewpoint,
$r$-th contact between manifolds was realised as the $r$-th contact of two
local parametrisations of the given manifolds at a point. Clearly, the bridge
from one viewpoint to another is provided by the notion of `graph': two maps
have an $r$-th contact if and only if their graphs have an $r$-th order
contact.

It turns out that the most general notion of contact between manifolds can be
given as contact between embedded submanifolds of a given manifold.  A
detailed justification of this statement is given in the subsection
\ref{ssec:jsubm_top}.

In the modern terminology, the spaces of equivalence classes of submanifolds of
a given manifold having an $r$-th contact at a certain point are said to be
\emph{spaces of jets of submanifolds} (also known as \emph{manifolds of contact
  elements}). The above arguments show that the geometry of jets of
submanifolds plays a key role whenever dealing with mathematical topics where
the geometry of contact is important, like geometric aspects of differential
equations and calculus of variations. Despite this fact, the preferred approach
to jets has been mostly through jets of maps or fibrings (see, \emph{e.g.},
\cite{Many99,KMS93,Kup80,MaMo83,Sau89}), because the last ones generalise the
notion of graph through the concept of `section'.  Instead, it is clear that
jets of submanifolds are locally diffeomorphic to jets of fibrings, so jets of
submanifolds constitute a non-trivial generalisation of jets of fibrings.

Historically, less attention has been devoted to jets of submanifolds with
respect to other theories of jets.  After the work of Lie, the research by
Bompiani\footnote{He used the name \emph{elementi differenziabili} to indicate
  jets of submanifolds} was devoted to the geometry of jets of surfaces about
1910 and during all his career \cite{Bom1}.  He was especially concerned in
projective invariants of jets of curves and surfaces and applications to
differential equations. It seems that his approach has been left apart. He also
wrote a remarkable work on the history of the notion of contact, which he dates
back to Ruffini~\cite{Bom2}.

Jets of submanifolds have been first introduced in a modern setting by C.
Ehresmann~\cite{Ehr52}. He defined the notion of contact of order $r$ between
parametrised submanifolds, then identified $r$-equivalent parametrisations. The
works~\cite{CrSa03a,CrSa03b,CrSa03c,GrKr98,Gug79,KMS93,Kol98,Kru01,MMR00}, are
among the few researches performed within this scheme. We also mention the
brief introduction in \cite{Olv91} (under the name of \emph{extended jet
  bundles}).  In the classical calculus of variations, jets of submanifolds
arise implicitly when considering parametric
problems~\cite{GiHi96}\footnote{Such problems date back to Carath\'eodory in
  the case of one independent variable}. In such variational problems the
Lagrangian is invariant with respect to changes of parametrisation, hence it
factors to the space of jets of submanifolds with respect to the action of
change of parametrisation.

More directly, in~\cite{Ded77,Vin77,Vin78,Vin84} jets of submanifolds are
introduced through the notion of contact of order $r$ between submanifolds. In
this paper, we followed this second approach. Indeed, we find this viewpoint
more straightforward, especially for what concerns coordinate expressions,
because the dependency on parametrisation is factored out. Authors working in
the parametrised scheme like~\cite{Gri97,Gri98,GrKr98,CrSa03a,CrSa03b,CrSa03c}
are forced to check each time the independence of their computations from
parametrisations.

\medskip

Our guideline in the study of jets of submanifolds has been to following:
\emph{to proceed by analogy with other theories of jets}, in particular jets of
fibrings, and try to reproduce similar geometric structures and constructions.

This is our second paper about jets of submanifolds.  In the previous
one~\cite{MaVi02} we introduced the pseudo-horizontal and pseudo-vertical
bundles on jets of submanifolds (see subsection~\ref{ssec:contactseq}). Such
bundles play the role of horizontal and vertical bundles in the case of jets of
fibrings. By those bundles we computed the cohomology of the finite order
analogue of Vinogradov's $\cC$-spectral sequence (see section~\ref{sec:C-ss}).

In this paper we begin with a more careful exposition of basic definitions
(subsections~\ref{ssec:jets}, \ref{ssec:contactseq}), together with a thorough
comparison between the three known approaches to jet spaces
(subsection~\ref{ssec:jsubm_top}). Then, we derive a new proof of
the fact that higher-order jets of submanifolds are affine bundles
(subsection~\ref{ssec:affine}). Other proofs in
literature~\cite{AVL91,Kol98,MMR00} rely on the analysis of the transformation
group of the fibres. Our proof is based on the structure of the
pseudo-horizontal and pseudo-vertical bundles, and is a generalisation of the
similar result in the case of jets of fibrings.  As a by-product, we are able
to provide a new explicit expression of the associated vector bundles using
tensor products of pseudo-horizontal and pseudo-vertical bundles.

Next, we recall briefly the computations that led in~\cite{MaVi02} to the
computation of the finite order variational sequence on jets of submanifolds.
The sequence contains spaces of forms with a special polynomial structure. Such
a structure has been studied in coordinates so far
(\emph{hyperjacobians},~\cite{Olv83}, and~\cite{AnDu80}, \cite[chap. 4]{And}).
In this paper we find it as a by-product of our geometric structure
(subsection~\ref{ssec:forms} and section~\ref{sec:C-ss}). We also provide a
useful coordinate expression for the first differential of the spectral
sequence, $\hd$ (proposition~\ref{pro:coordex_hd}).

Then, we show that a finite order version of Green-Vinogradov formula does not
exist (subsection~\ref{ssec:green-vinogr-form}). This leads us to use the
infinite-order theory to compute the representative of each cohomology class
(also known as \emph{variational form}) in the variational sequence. Such
representatives are computed in theorem~\ref{th:iso} through~\eqref{eq:repr}.
As a by-product, we solve a problem from~\cite{Gra00} (conjecture/question by
P. Griffiths, remark~\ref{re:hor_coh}). This is due to the fact that we do not
confine ourselves to the computation of the horizontal cohomology (\emph{i.e.},
the cohomology of the horizontal de Rham
sequence~\cite{Many99,KLV86,KrVe98,Ver98}), but we compute all cohomology
groups of the finite order $\cC$-spectral sequence.

We provide the coordinate expression for all representatives of variational
forms in subsection~\ref{ssec:coord-expr}. In the cases $p=0,1,2$ these
expressions look like the well-known expressions of Lagrangians,
Euler--Lagrange operators and Helmholtz operators of the case of jets of
fibrings, but the spaces where they are defined are quite different, hence also
their coordinate transformation laws. As far as we know, it is the first time
that such expressions are computed on jets of submanifolds.

Finally, in section~\ref{sec:var_princ} we provide a geometric (\emph{i.e.},
invariant, coordinate-free) formulation of the variational problem on jets of
submanifolds. This allows us to derive the Euler--Lagrange equations in an
intrinsic way. Our formulation reduces to well-known formulations in the case
of jets of fibrings (see~\cite{KMS93,Kru01,Kup80,Sau89,Tul77}, for example),
and is a radical improvement of the old formulation by Dedecker~\cite{Ded77}.
We also describe the connection of this formalism with the finite order
$\cC$-spectral sequence. Indeed, we think that this formalism provides one of
the main motivations for the $\cC$-spectral sequence itself.

We already observed that variational problems on jets of submanifolds have been
introduced in the parametric
framework~\cite{CrSa03a,CrSa03b,CrSa03c,Gri97,Gri98}. In our opinion, we
provided a simpler formulation through pseudo-horizontal bundles. As an
example, we deal with a single Lagrangian, not with a set of Lepage
equivalents. Moreover, we do not need to check the invariance with respect to
changes of parametrisations.  We added detailed comparisons with literature in
section~\ref{sec:var_princ} explaining our viewpoint.

As a final remark, we regret that we had no time and space to support our
theory with examples, like in~\cite{CrSa03a}. We have some partial
results~\cite{Man03}, that we hope to work out and expose in a next paper.
\section{Jet spaces}\label{sec:jets}

In this section we recall basic facts about the geometry of jets of
submanifolds (our sources were~\cite{AVL91,MoVi94,Vin88}) together with new
structures (subsections~\ref{ssec:contactseq}, \ref{ssec:affine}) which
will be used for the analysis of the finite order $\cC$-spectral sequence.

\medskip

In this paper all manifolds and maps are smooth, and all submanifolds are
\emph{embedded submanifolds}.

We shall introduce jets through the fundamental notion
of contact between maps. This is sometimes given through composition
of maps with curves (see, \emph{e.g.},~\cite{KMS93}); here, we adopt a more
direct viewpoint.

Let $N$ and $M$ be two manifolds, and $f,g\colon N\to M$ two maps.  We
say that $f$ and $g$ have a \emph{contact of order $r$} at $x$ if
there is a chart $(U,(x^\lambda))$ at $x$ and a chart $(V,(u^i))$ at
$f(x)$ such that
\[
J^k((u^i)\circ f\circ(x^\lambda)^{-1})=J^k((u^i)\circ
g\circ(x^\lambda)^{-1}),
\qquad 0\leq k\leq r,
\]
where $J^k$ stands for the $k$-th order differential. Of course,
$f$ and $g$ have a contact of order zero if $f(x)=g(x)$.

It is easy to realize that the above notion is independent of the
chosen chart.
\subsection{Jet spaces}\label{ssec:jets}

Let $E$ be an $(n+m)$-dimensional manifold.

Let $(V,(y^h))$ be a local chart on $E$.  Any splitting
$\{1,\ldots,n+m\}=I_n\cup I_m$ into two subsets $I_n$ and $I_m$
having, respectively, $n$ and $m$ elements, induces a splitting
$(y^h)=(x^\lambda,u^i)$ of the chart, where
$x^\lambda=y^{h_{\lambda}}$, $u^i=y^{h_{i}}$ and $h_\lambda\in
I_n$, $h_i\in I_m$, so that $1\leq \lambda\leq n$ and $1\leq i\leq
m$.  The above splitting is said to be a \emph{division} of the
chart $(y^h)$. A chart of the form $(x^\lambda,u^i)$ is said to be
a \emph{divided
  chart}. Of course, there are $\binom{n+m}{n}$ possible divisions of $(y^h)$.

In what follows, Greek indices run from $1$ to $n$ and Latin indices run from
$1$ to $m$.

We say that a divided chart $(V,(x^\lambda,u^i))$ is \emph{fibred} if $V$ is
diffeomorphic to $X \times U$, where $X\subset{\R^n}$ and $U\subset{\R^m}$ are
open subsets. Of course, the trivial projection $\pi\colon V\to X$ makes $V$ a
fibred manifold on $X$.

Let $L$ be an (embedded) submanifold $\iota\colon L \hookrightarrow E$ ($\iota$
is the inclusion map). We say that a divided chart $(V,(x^\lambda,u^i))$ is
\emph{concordant} with $L$ at $x\in V\cap L$ if the coordinate expression of
$\iota$ is
\[
(x^\lambda,u^i)\circ\iota=(x^\lambda,\iota^i),
\]
where the functions $\iota^i$ are smooth functions of
$(x^\lambda)$. We say that the chart is \emph{adapted} to $L$ at
$x$ if $(\iota^i)=(0)$. Such a chart exists at any $x\in L$.

Let $\iota\colon L\hookrightarrow E$, $\iota'\colon
L'\hookrightarrow E$ be two submanifolds, and $x\in L\cap L'$. Let
$\bsi =(\sigma _{1},\sigma _{2},\ldots,\sigma_{r})$, with
$1\leq\sigma_i\leq n$ and $r\in\mathbb{N}$, be a multi-index, and
$\left| \bsi \right|\byd r$.  Then, we say that $L$ and $L'$ have
a \emph{contact of order $r$} at $x$ if $\iota$ and $\iota'$ have
a contact of order $r$ at $x$.

Let us examine the above definition in coordinates. It can be proved that there
exists a fibred chart which is concordant to both $L$ and $L'$ (by the theory
of transversality, see, \emph{e.\ g.},~\cite{Hir76}).  If $(x^\lambda,u^i)$ is
such a chart at $x$, then $L$ and $L'$ have a contact of order $r$ at $x\in
L\cap L'$ if
\begin{equation}\label{eq:jet}
\pd{^{\abs{\bsi}}\iota^i}{x^{\sigma_1}\cdots\partial x^{\sigma_r}}(x)=
\pd{^{\abs{\bsi}}{\iota'}^i}{x^{\sigma_1}\cdots\partial
  x^{\sigma_r}}(x),
\qquad 0\leq\abs{\bsi}\leq r.
\end{equation}
The relation ``contact of order $r$ between submanifolds at $x\in E$'' is an
equivalence relation; an equivalence class is denoted by~$[L]^{r}_{x}$.
\begin{definition}
  Let
  \begin{displaymath}
    J^r(E,n)\byd\bigcup_{x\in E} J^r_x(E,n),
  \end{displaymath}
  where $J^r_x(E,n)$ is the set of the equivalence classes $[L]^r_x$ of
  $n$-dimensional submanifolds $L\subset E$ having a contact of order $r$ at
  $x$. We call $J^r(E,n)$ the \emph{$r$-jet of $n$-dimensional
    submanifolds of $E$}.
\end{definition}
Note that $J^0(E,n)=E$. Any submanifold $L \subset E$ can be prolonged
to a subset of $J^r(E,n)$ through the injective map
\begin{displaymath}
j_{r}L\colon L \to J^{r}(E,n),\quad x \mto [L]_{x}^{r}.
\end{displaymath}
The set $J^r(E,n)$ has a natural manifold structure: any divided
chart $(V,x^\lambda,u^i)$ induces the local chart
$\left(V^r_n,(x^{\lambda },u_{\bsi}^{i})\right)$ on $J^{r}(E,n)$
at $[L]_{x}^{r}$, where $V^r_n\byd\bigcup_{x\in V} J^r_x(E,n)$,
and the functions $u^j_{\bsi}$ are determined by
\begin{equation*}
u_{\bsi}^{i}\circ j_{r}L=\pd{^{\left| \bsi \right| }\iota^{i}%
}{x^{\sigma_1}\cdots\partial x^{\sigma_r}}.
\end{equation*}
Coordinate changes are smooth because they involve compositions
and sums of $k$-th order Jacobians ($0\leq k\leq r$) of the
coordinate change in $E$.  The dimension of $J^r(E,n)$ is the sum
of the dimension of $E$ with the number of all the possible
derivatives of $u^i$ with respect to $x^\lambda$ up to order $r$,
hence
$\dim J^r(E,n)=n+m\sum_{h=0}^r\binom{n+h-1}{h}=n+m\binom{n+r}{r}$.

The map $j_rL$ turns out to be an embedding; we shall
identify it with its image, denoting it by $L^{(r)}$.

Any local diffeomorphism $F\colon E\to E$ prolongs to a fibred isomorphism
$J^rF\colon$ $J^r(E,n)$ $\to$ $J^r(E,n)$ in a way similar to jets of fibrings
(see, for example,~\cite{MaMo83,Sau89}). Namely, we have
$J^rF([L]^r_x)=[F(L)]^r_{F(x)}$. Hence, any vector field $X\colon E\to TE$
prolongs to a vector field $X_r\colon J^r(E,n)\to TJ^r(E,n)$, by prolonging its
flow.

The above smooth manifold structure endows the natural projections
\begin{equation*}
\pi_{r,h}:J^{r}(E,n)\to J^{h}(E,n),\quad r\geq h,
\end{equation*}
with a bundle structure. In particular, it is
known~\cite{AVL91,Kol98,Vin84} that $\pi_{r+1,r}$ are affine
bundles for $r\geq 1$. Later on, we shall provide a more complete
proof in our framework. We shall also denote by $J^\infty (E,n)$
the inverse limit of the chain of projections
\[
\begin{diagram}
  \cdots & \rTo^{\pi_{r+1,r}} & J^r(E,n) & \rTo^{\pi_{r,r-1}} & \cdots &
  \rTo^{\pi_{2,1}} & J^1(E,n) & \rTo^{\pi_{1,0}} & E.
\end{diagram}
\]

A special attention is needed in the case $r=0$. It is easy to realize that
$\pi_{1,0}$ coincides with the Grassmann bundle of $n$-dimensional subspaces in
$TE$. Any chart at $x\in E$ induces a covering of the Grassmann space $\Gr
(T_xE,n)\simeq \pi_{r,0}^{-1}(x)$ $(x^\lambda,u^i)$ as follows. Any division
$(x^\lambda,u^i)$ of the given chart induces a local chart of $\Gr
(T_xE,n)$ onto the open set made by $n$-dimensional
subspaces spanned by the vectors
\begin{equation}\label{eq:Grass}
\eval{\pd{}{x^\lambda}}_x+u^i_\lambda\eval{\pd{}{u^i}}_x.
\end{equation}
Each of the above open subsets is dense in $\Gr (T_xE,n)$, and the set of all
divisions of the given chart covers $\Gr (T_xE,n)$

Note that $\pi_{1,0}$ has no affine structure (this contrasts with the case of
jets of fibrings, see~\cite{MaMo83,Sau89}). It is easy to realize that
$\pi_{1,0}$ has a Grassmannian structure, in the sense that fibred coordinate
changes are Grassmannian transformations induced by isomorphisms of $TE$. More
precisely, let $(x^\lambda,u^i)$ and $(y^\mu,v^j)$ be two coordinate charts
concordant to the same submanifold $L\subset E$; let us denote by
$(J_\lambda^\mu,J_i^\mu,J_\lambda^j,J_i^j)$ the Jacobian of the change of
coordinates. Then the fibred coordinate change is given by the following
formula
\begin{equation}\label{eq:changefirst}
v^j_\mu=\frac{J^j_\lambda+J^j_i u^i_\lambda}
{J^\mu_\lambda+J^\mu_i u^i_\lambda}.
\end{equation}
In fact, locally, the submanifold $L$ is expressed by
$u^i=f^i(x^\lambda)$ or $v^j=g^j(y^\mu)$. So, differentiating the
equation $v^j|_L=v^j(x^\lambda,f^i(x^\lambda))=
g^j\left(y^\mu(x^\lambda,f^i(x^\lambda))\right)$ we get the
result.
\subsection{Jets of submanifolds as a generalisation of jets of fibrings
 and maps}\label{ssec:jsubm_top}

We hereafter provide justifications to the title of this
subsection, of intrinsic and local nature. Suppose that $E$ is
endowed with a fibring $\pi :E\to M$.
\begin{itemize}
\item The space $J^r\pi$ of $r$-th jets of sections $s\colon M\to E$ of $\pi$
  is an open dense subset of $J^r(E,n)$. In fact, it coincides with the subset
  of $r$-th jets of submanifolds of the type $s(M)$, and it is covered by just
  one of the open subsets of the previous covering of $\Gr (T_xE,n)$
  \eqref{eq:Grass}. See~\cite{Many99,KrVe98,MaMo83,Sau89} for the theory of
  jets of fibrings.
\item When $E=X\times U$, then there is the trivial fibring $\pro\colon X\times
  U\to X$, and $J^r\pro$ coincides with the space of $r$-jets of maps
  $J^r(X,U)$. In fact, the space $J^r(X,U)$ is defined as the set of
  equivalence classes of maps of the type $f\colon X\to U$ having an
  $r$-contact at a point, hence we fall back in the previous case.
  See~\cite{Hir76,KMS93} for definitions and properties of jets of maps.
\item A comparison between transformation rules of independent and dependent
  variables in various approaches to jets also shows that jets of submanifolds
  are the most general theory~\cite{Ded77,Vin88}. Let $(x^\lambda,u^i)$ and
  $(\bar{x}^\lambda,\bar{u}^i)$ be two charts at the same point. Then we have
  the transformation rules:
\begin{itemize}
\item $(x^\lambda(\bar{x}^\mu), u^i(\bar{u}^j))$ (jets of maps),
\item $(x^\lambda(\bar{x}^\mu), u^i(\bar{x}^\mu,\bar{u}^j))$ (jets of
  fibrings),
\item $(x^\lambda(\bar{x}^\mu,\bar{u}^j), u^i(\bar{x}^\mu,\bar{u}^j))$ (jets of
  submanifolds).
\end{itemize}
\end{itemize}

It should be said that $J^r(E,n)$ can be derived from jets of maps as follows.
Consider the subspace $\Imm J^r_0(\R^n,E)\subset J^r(\R^n,E)$ of $r$-jets of
immersions at $0\in\R^n$.  Any immersion can be seen as a local parametrisation
of an (embedded) submanifold of $M$. There is an action of the group $G^r_n$ of
$r$-jets of local diffeomorphisms $f$ of $\R^n$ such that $f(0)=0$ on $\Imm
J^r_0(\R^n,E)$; this can be regarded as the action by changes of
parametrisation. Then it is easy to see that the following isomorphism holds
\[
\Imm J^r_0(\R^n,E)/G^r_n\simeq J^r(E,n).
\]
This can be taken as a definition of jets of submanifolds: see
\cite{Gri97,Gri98,GrKr98,KMS93,Kru01}.
\begin{remark}
  Some authors prefer to work on $\Imm J^r_0(\R^n,E)$ with objects that factor
  through the action by $G^r_n$ to $J^r(E,n)$
  \cite{Gri97,Gri98,GrKr98,CrSa03a,CrSa03b,CrSa03c}. This amounts to work with
  $n$ extra parameters and to show each time the independence of the
  parameters. In our opinion, this method seems to be more involved with
  respect to the direct use of jets of submanifolds. We especially
  get results in a more straightforward way with respect to the
  `parametric' approach in section~\ref{sec:var_princ}.
\end{remark}
\subsection{Contact sequence}
\label{ssec:contactseq}

Here, we introduce a new definition of contact structure on jets which
turns out to be very useful in some cases, like in the analysis of finite order
variational sequences~\cite{MaVi02}. This idea generalises the definition that
is given in~\cite{MoVi94} for first-order jets, but it is somehow
different from the standard contact (or Cartan) distribution on jets
which is of fundamental importance for the geometric study of
differential equations~\cite{Many99}.

\bigskip

For $r\geq 0$ consider the following bundles over $J^{r+1}(E,n)$: the
pull-back bundle
\begin{equation}
T^{r+1,r}\byd J^{r+1}(E,n)\underset{J^{r}(E,n)}{\times}T
J^{r}(E,n),
\end{equation}
the subbundle $H^{r+1,r}$ of $T^{r+1,r}$ defined by
\begin{equation}
H^{r+1,r}\byd\left\{ \left( [L]_p^{r+1} ,\upsilon
\right) \in T^{r+1,r} \mid \upsilon \in T_{[L]_p^r}L^{\left( r\right) }\right\}
\end{equation},
and the quotient bundle
\begin{equation}
V^{r+1,r}\byd T^{r+1,r}/H^{r+1,r}.
\end{equation}

The bundles $H^{r+1,r}$ and $V^{r+1,r}$ are strictly related with the
horizontal and vertical bundle in the case of jets of a fibring (see
remark~\ref{re:fibr}).

\begin{definition}
We call $H^{r+1,r}$ and $V^{r+1,r}$, respectively, the
\emph{pseudo-horizontal} and the \emph{pseudo-vertical} bundle of
$J^r(E,n)$.
\end{definition}

First of all, we see how the above bundles relate to the projections
$\pi_{r+1,r}$. There is a natural projection $T^{r+1,r}\to T^{r,r-1}$.
It restricts to $H^{r+1,r}\to H^{r,r-1}$, hence it induces a natural
projection $V^{r+1,r}\to V^{r,r-1}$.

We observe that, while $(T^{r+1,r})^*$ is a vector subbundle of
$T^*J^{r+1}(E,n)$ via the inclusion $T^*\pi_{r+1,r}$, there is no
natural inclusion of $T^{r+1,r}$ into $TJ^{r+1}(E,n)$.

The pseudo-horizontal bundle has some additional features. The
following isomorphism over $\id_{J^{r+1}(E,n)}$ holds
\begin{equation}\label{eq:1}
H^{r+1,r} \to J^{r+1}(E,n)\times _{J^{1}(E,n)}H^{1,0},\qquad
([L] _{x}^{r+1},\upsilon)\mto
\left([L] _{x}^{r+1},T\pi _{r,0}(\upsilon )\right).
\end{equation}
Hence, we obtain the natural projection $(H^{r+1,r})^*\to(H^{r,r-1})^*$.

Note that, if $L\subset E$ is an $n$-dimensional submanifold of $E$,
then the pull-back bundle $(j_{r+1}L)^*H^{r+1,r}$ on $L$ is just $TL$,
as it is easily seen. This is a further justification for its name:
indeed, the pseudo-horizontal bundle is tangent to the jet of any
$n$-dimensional submanifold of $E$.

The most important property of the bundles $T^{r+1,r}$,
$H^{r+1,r}$ and $V^{r+1,r}$ is the following \emph{contact exact
  sequence}
\begin{equation}\label{eq:contact}
\newdiagramgrid{contact}{1,.8,1,1.3,1,1.3,1,.8,1}{1}
\begin{diagram}[grid=contact]
0 & \rTo & H^{r+1,r} & \rhTo^{D^{r+1}} & T^{r+1,r} &
\rTo^{\omega^{r+1}} & V^{r+1,r} & \rTo & 0,
\end{diagram}
\end{equation}
where $D^{r+1}$ and $\omega^{r+1}$ are the natural inclusion and
quotient projection. It induces the following exact sequence
\begin{equation}
  \label{eq:2}
\newdiagramgrid{contseq}
{1,.6,1.2,2,2.2,2.1,2,.5,1}{1}
\begin{diagram}[grid=contseq]
  0 & \lTo & \bigwedge^k(H^{r+1,r})^* & \lTo^{\bigwedge^k(D^{r+1})^*} &
  \bigwedge^k(T^{r+1,r})^* & \lhTo^{(\omega^{r+1})^*\wedge\,\id} &
  (V^{r+1,r})^*\wedge\bigwedge^{k-1}(T^{r+1,r})^* & \lTo & 0.
\end{diagram}
\end{equation}

\begin{remark}\label{re:fibr}
  When $E$ is endowed with a fibring $\pi :E\to M$, we have the following
  isomorphisms over $J^{r+1}\pi$:
\[
\begin{diagram}
H^{r+1,r}|_{J^{r+1}\pi} & \rTo & J^{r+1}\pi\times_M TM,
\end{diagram}\qquad
\begin{diagram}
V^{r+1,r}|_{J^{r+1}\pi} & \rTo & J^{r+1}\pi\times_{J^r\pi}VJ^r\pi,
\end{diagram}
\]
where $VJ^r\pi\byd\ker T\pi_{r,0}$. Hence the exact sequence over
$J^{r+1}\pi$
\[
\newdiagramgrid{horseq}{1,1,2,1.5,2,1.7,2,.3,1}{1}
\begin{diagram}[grid=horseq]
0 & \rTo & J^{r+1}\pi\times_{J^r\pi}VJ^r\pi & \rhTo &
J^{r+1}\pi\times_{J^r\pi} TJ^r\pi & \rTo^{T(\pi\circ\pi_{r,0})} &
J^{r+1}\pi\times_M TM & \rTo & 0
\end{diagram}
\]
splits the contact sequence~\eqref{eq:contact} \cite{MaMo83,Sau89}.
\end{remark}

Let us evaluate the coordinate expressions of $D^{r+1}$ and $\omega^{r+1}$.
 To this end, we observe that
$D^{r+1}$ and $\omega ^{r+1}$ can be seen as sections of the
bundles $\left(H^{r+1,r}\right)^{\ast }\otimes_{J^{r+1}(E,n)}
T^{r+1,r}$ and $\left(T^{r+1,r}\right)^*\otimes_{J^{r+1}(E,n)}
V^{r+1,r}$ respectively.

A local basis of the space of sections of the bundle $H^{r+1,r}$ is
\begin{equation*}
D^{r+1}_\lambda=\frac{\partial }{\partial x^{\lambda }}+
u_{\bsi,\lambda}^{j}\frac{\partial }{\partial u_{\bsi}^{j}},
\end{equation*}
where the index $\bsi,\lambda$ stands for
$(\sigma_1,\dots,\sigma_s,\lambda)$ with $s\leq r$.
A local basis of the space of sections $(H^{r+1,r})^*$ dual to
$(D^{r+1}_\lambda)$ is given by the restriction of the $1$-forms $dx^\lambda$
to $H^{r+1,r}$, and is denoted by $\overline{dx}{}^\lambda$.
The local expression of $D^{r+1}$ turns out to be
\begin{equation*}
D^{r+1}=\overline{dx}{}^\lambda\otimes
D^{r+1}_\lambda=\overline{dx}{}^\lambda\otimes
 \left( \frac{\partial }{\partial x^{\lambda }}%
+u_{\bsi,\lambda}^{j}\frac{\partial }{\partial u_{\bsi}^{j}}%
\right).
\end{equation*}
A local basis of the space of sections of the bundle $V^{r+1,r}$ is
\[
B^{\bsi}_{j}\byd\left[\frac{\partial}{\partial u^j_{\bsi}}
\right]\,\,, \qquad \abs{\bsi}\leq r.
\]
The local expression of $\omega^{r+1}$ is
\begin{equation*}
\omega ^{r+1}=\omega^j_{\bsi}\otimes B^{\bsi}_{j}=
\left( du_{\bsi}^{j}-u_{\bsi,\lambda}^{j}dx^{\lambda
}\right) \otimes B^{\bsi}_{j}.
\end{equation*}

\begin{remark}\label{re:Cartan}
  We have a natural distribution $\cC^r$ on $J^r(E,n)$ generated by the tangent
  spaces $TL^{(r)}$ for any $n$-dimensional submanifold $L\subset E$, namely
  the \emph{Cartan distribution}~\cite{AVL91,Many99}. It is generated by the
  vector fields $D^r_\lambda$ and $\partial/\partial u^i_{\bsi}$, with
  $\abs{\bsi}=r$.  This distribution has not to be confused with $H^{r,r-1}$,
  which is a subbundle of a different bundle and is generated by $D^r_\lambda$.
\end{remark}
\subsection{Bundle structures}
\label{ssec:affine}

Here we give a proof of the fact that $(J^{r+1}(E,n),\pi_{r+1,r},J^{r}(E,n))$
are affine bundles if $r\geq 1$. This is already known
\cite{AVL91,Kol98,MMR00}, but our proof is different and has the advantage to
provide also a new expression of the associated vector bundle in terms of our
pseudo-horizontal and pseudo-vertical bundles.

We stress that the proof is almost completely analogous to the case of jets of
fibrings \cite{MaMo83}, having introduced the analogues of horizontal and
vertical bundles in the previous subsection. The only difference is the absence
of the base space in our case, giving the obstruction to $\pi_{1,0}$ to be
affine.

\begin{lemma}\label{le:vert_iso}
The following isomorphism holds
\[
VJ^1(E,n)\simeq (H^{1,0})^*\otimes_{J^{1}(E,n)}V^{1,0}.
\]
\end{lemma}
\begin{proof}
  Any point of $J^1(E,n)$ can be seen as the inclusion of an $n$-dimensional
  subspace of $TE$ into $TE$ itself through $D^1$, hence as a linear map
  $\overline{dx}{}^\lambda\otimes (\pd{}{x^{\lambda}}+u^i_\lambda \pd{}{u^i})$.
  A curve tangent to the fibre of $\pi_{1,0}$ at such a point has the tangent
  vector $\overline{dx}{}^\lambda\otimes \dot u^i_\lambda \pd{}{u^i}$; this
  proves the above isomorphism.
\end{proof}
\begin{theorem}
  For $r\geq 1$ the bundles $(J^{r+1}(E,n),\pi_{r+1,r},J^{r}(E,n))$
  are affine bundles associated with the vector bundle
  \[
  \left(\left(\odot^{r+1}(H^{1,0})^*\right)
  \otimes_{J^r(E,n)}V^{1,0},\pro,J^r(E,n)\right),
  \]
  where $\pro$ is the trivial pull-back projection.
\end{theorem}
\begin{proof} We shall prove the theorem in two steps, dealing separately with
  the cases $r=1$ and $r>1$.

  For $r=1$ we can interpret $D^1$ as the section
  \[
    D^{1}\colon J^{1}(E,n)\to (H^{1,0})^*\otimes_{J^{1}(E,n)}TE.
  \]

  For $r>1$ we can interpret $D^{r}$ as a fibred inclusion
  \[
    D^{r}\colon J^r(E,n)\to (H^{1,0})^*\otimes_{J^{r-1}(E,n)}TJ^{r-1}(E,n)
  \]
  of bundles over $\id_{J^{r-1}(E,n)}$ through the isomorphism
  \eqref{eq:1}.

For $r=1$ we have the following commutative diagram
\begin{diagram}
  J^2(E,n) & \rTo^{D^2} & (H^{1,0}){}^*\otimes_{J^1(E,n)} TJ^1(E,n)
  \\
  \dTo^{\pi_{2,1}} & & \dTo_{\id\otimes T\pi_{1,0}}
  \\
  J^1(E,n) & \rTo^{D^1} & (H^{1,0}){}^*\otimes_{J^1(E,n)} TE
\end{diagram}
The inverse image through $\id\otimes T\pi_{1,0}$ of the
non-vanishing section $D^1(J^1(E,n))$ is an affine bundle. In view
of lemma \ref{le:vert_iso} the associated vector bundle is
\[
\ker(\id\otimes T\pi_{1,0})=(H^{1,0})^*\otimes_{J^1(E,n)}VJ^1(E,n)
\simeq
(H^{1,0})^*\otimes_{J^1(E,n)}(H^{1,0})^*\otimes_{J^1(E,n)}V^{1,0}.
\]
The above affine bundle contains $J^2(E,n)$ as the subbundle
$D^2(J^2(E,n))$. It is easy to realise that $J^2(E,n)\to J^1(E,n)$
is the affine subbundle whose associated vector bundle is
\[
\left((H^{1,0})^*\odot_{J^1(E,n)}(H^{1,0})^*\right)
\otimes_{J^1(E,n)}VJ^1(E,n).
\]

For $r>1$ we have the following commutative
diagram
\begin{diagram}
  J^{r+1}(E,n) & \rTo^{D^{r+1}} & (H^{1,0})^*\otimes_{J^{r}(E,n)}
  TJ^{r}(E,n)
  \\
  \dTo^{\pi_{r+1,r}} & & \dTo_{\id\otimes T\pi_{r,r-1}}
  \\
  J^{r}(E,n) & \rTo^{D^{r}} & (H^{1,0})^*\otimes_{J^r(E,n)}
  TJ^{r-1}(E,n)
\end{diagram}
and we get that $\pi_{r+1,r}$ is an affine bundle by a
similar reasoning. It is now easy to obtain the associated vector bundle
by induction.
\end{proof}
\subsection{Forms on jets}
\label{ssec:forms}

Here we study the spaces of forms on jets in view of the
formulation of the finite order $\mathcal{C}$-spectral sequence.
We introduce spaces of \emph{contact} and \emph{horizontal} forms.
Contact forms vanish when calculated on any prolonged submanifold;
horizontal forms are forms which do not contain contact factors.
We show that horizontal forms have a special polynomial structure
which has already studied from an intrinsic viewpoint in the case
of jets of fibrings~\cite{PaVi00,Vit98}. We prove that such a
structure is present also on horizontal forms on zero-order jets.

Note that such a polynomial structure has been introduced and studied in a
coordinate (or local) fashion so far (\emph{hyperjacobians},~\cite{Olv83},
and~\cite{AnDu80}, \cite[chap. 4]{And}); in this paper we find it as a
by-product of our geometric structure. This also provides a better
understanding of its transformation laws.

\medskip

We denote by $\cF_r$ the algebra $C^{\infty}(J^r(E,n))$.  For
$k\geq 0$ we denote by $\Lambda^{k}_r$ the $\cF_r$-module  of
$k$-forms on $J^r(E,n)$. We also set $\Lambda^{*}_r=\bigoplus
_{k}\Lambda^{k}_r$.  We introduce the submodule of
$\Lambda_{r}^{k}$ of the \emph{contact forms}
\begin{equation*}
\cC^{1}\Lambda _{r}^{k}\byd\{\,\alpha \in \Lambda _{r}^{k} \mid
(j_{r}L)^{\ast }\alpha =0\quad\text{for each submanifold
}L\subset E\,\}.
\end{equation*}
Contact forms are clearly the annihilators of the Cartan distribution (see
remark~\ref{re:Cartan}).  We set $\cC^1\Lambda^*_r=\bigoplus_k
\cC^1\Lambda^k_r$. Moreover, we define $\cC^p\Lambda^*_r$ as the
$p$-th exterior power of $\cC^1\Lambda^*_r$. Of course,
$\cC^p\Lambda^{k}_r=\cC^p\Lambda^*_r\cap\Lambda^k_r$.

Next we introduce the $\cF_{r+1}$-module $\Lambda^k_{r+1,r}$ of
sections of the bundle $\bigwedge^k(T^{r+1,r})^*$. This module is formed
by $k$-forms along $\pi_{r+1,r}$, \emph{i.e.}, $k$-forms on $J^r(E,n)$ with
coefficients in $\cF_{r+1}$.

We also consider the $\cF_{r+1}$-module $\cH^k_{r+1,r}$ of
\emph{pseudo-horizontal} $k$-forms, \emph{i.e.}, sections of the
bundle $\bigwedge^k(H^{r+1,r})^*$.

\medskip

Now we introduce an operation that allows us to extract from any form
$\alpha\in \Lambda_{r}^{k}$ its \lq horizontal part\rq.

\begin{definition}\label{def:horizon}
Let $q\in\mathbb{N}$. \emph{Horizontalisation} is the map
\begin{equation*}
  h^{0,q}\colon   \Lambda _{r}^{q}  \to \cH_{r+1,r}^{q}, \quad
  \alpha \mto (\wedge^q(D^{r+1})^*)\circ (\pi _{r+1,r}^*\alpha).
\end{equation*}
\end{definition}
Of course, the above operation is trivial if $q>n$.
Horizontalisation is well-defined because
$\pi_{r+1,r}^*(\alpha)$ has values in $\bigwedge^q(T^{r+1,r})^*$,
which is a subbundle of $\bigwedge^q(T^*J^{r+1}(E,n))$ (see
subsection~\ref{ssec:contactseq}). If $\alpha\in\Lambda^{1}_{r}$
has the coordinate expression $\alpha =\alpha_\lambda dx^\lambda
+\alpha_i^{\bsi}du^i_{\bsi}$ ($0 \leq |{\bsi}| \leq r$), then
\begin{displaymath}
h^{0,1}(\alpha ) =(\alpha_\lambda+u^i_{\bsi,\lambda}\alpha_i^{\bsi})\
\overline{dx}{}^\lambda.
\end{displaymath}
In general if $\alpha\in\Lambda^q_r$, then we have the coordinate
expression
\begin{equation}\label{eq:3}
\alpha = \alpha
{_{i_1 \dots i_{h} }^{\bsi_1 \dots \bsi_{h}}}
{_{\lambda_{h+1} \dots \lambda_{q}}}
du^{i_1}_{\bsi_1}\wedge\dots\wedge du^{i_{h}}_{\bsi_{h}}\wedge
dx^{\lambda_{h+1}} \wedge\dots\wedge dx^{\lambda_{q}},
\end{equation}
where $0 \leq h \leq q$. Hence
\begin{equation}\label{eq:4}
h^{0,q}(\alpha) = u^{i_1}_{\bsi_{1},\lambda_{1}} \dots
u^{i_{h}}_{\bsi_{h},\lambda_{h}} \alpha
{_{i_1 \dots i_{h} }^{\bsi_1 \dots \bsi_{h}}}
{_{\lambda_{h+1} \dots \lambda_{q}}}
\overline{dx}{}^{\lambda_{1}} \wedge\dots\wedge \overline{dx}{}^{\lambda_{q}}.
\end{equation}

Let us introduce the $\cF_r$-module $\overline{\Lambda }_{r}^{q}\byd\im
h^{0,q}$. It is easy to realize from the above coordinate expressions that, if
$r\geq 1$, then $\overline{\Lambda }_{r}^{q}$ is made by elements of
$\cH^q_{r+1,r}$ whose coefficients are fibred polynomials of degree $q$ in the
highest order variables (i.\ e., $u^i_{\boldsymbol{\sigma}}$ with
$\abs{\boldsymbol{\sigma}}=r+1$).  Of course, this feature is intrinsic due to
the affine structure of $\pi_{r+1,r}$. But $\overline{\Lambda}_{r}^{q}$ does
not coincide with the space of \emph{all} such polynomial forms: indeed, not
all polynomial forms come from the horizontalisation of a form on a jet space,
unless $n=1$. From \eqref{eq:4} it follows that the coefficients of monomials
present a skew-symmetry with respect to the exchange of pairs ${}^i_{\bsi}$ and
${}^j_{\bta}$. This property appears analogously in the case of jets of
fibrings, see~\cite{Vit98,Vit99,Vit01}.

The case $r=0$ needs a special attention. In a similar way to
\eqref{eq:changefirst}, we realize that $\overline{dy}^\mu =
(J^\mu_\lambda+J^\mu_i u^i_\lambda)\overline{dx}^\lambda$.
Combining this formula with \eqref{eq:4} we deduce that the set of
sections of the bundle $\overline{\Lambda}_{0}^{q}$ admit a
subspace of sections with polynomial coefficients.  So, even if
$\pi_{1,0}$ \emph{is not} an affine bundle,
$\overline{\Lambda}_{0}^{q}$ is a subspace of the space of forms
with polynomial coefficients of degree $q$ in $u^i_\lambda$.
\begin{remark}\label{re:hyperjac}
  The above polynomial structure has been studied in a coordinate setting
  in~\cite{Olv83}. The alternated sum of monomials of the type
  $u^{i_1}_{\bsi_{1},\lambda_{1}} \dots u^{i_{h}}_{\bsi_{h},\lambda_{h}}$ is
  called \emph{hyperjacobian}.  We stress, however, that such a structure
  emerges naturally by virtue of the geometric properties of our scheme, hence
  it is a \emph{global} property, in contrast with the analysis
of~\cite{AnDu80} and~\cite[chap. 4]{And}, where such a polynomial structure is
  treated as a local property.
\end{remark}

Then, we study contact forms and their relationship with
horizontalisation.
\begin{lemma}\label{lem:app}
  Let $\alpha\in\Lambda^q_r$, with $0\leq q\leq n$. We have
  $(j_{r}L)^{\ast}(\alpha)=(j_{r+1}L)^{\ast}(h^{0,q}(\alpha ))$.
\end{lemma}
\begin{proof}
  This follows from the expression of $h^{0,q}$ (definition~\ref{def:horizon})
  and the fact that $D^{r+1}$ is the identity on the image of $Tj_rL$,
  \emph{i.e.}, on $TL^{(r)}$.
\end{proof}
As an obvious consequence of the previous lemma, we have
\begin{equation}\label{cC_and_ker_h}
\cC^1\Lambda^q_{r} = \ker h^{0,q} \quad \text{if} \quad 0\leq q \leq
n, \qquad \cC^1\Lambda^q_{r} = \Lambda^q_{r} \quad \text{if} \quad
q > n.
\end{equation}
Moreover, the forms in $\cC^1\Lambda^q_{r}$ can be characterised as follows:
\begin{equation}
\alpha \in \cC^1\Lambda^q_{r}
\quad\Leftrightarrow\quad
\pi_{r+1,r}^*(\alpha)\in\im
  ((\omega^{r+1})^*\wedge\,\id).
\end{equation}
It turns out that, if $\alpha \in \cC^p\Lambda^{p+q}_{r}$, then
we have the coordinate expression
\begin{equation}\label{eq:coord_contact}
\pi^*_{r+1,r}(\alpha)=\omega^{i_1}_{\bsi_1}\wedge\dots
\wedge\omega^{i_p}_{\bsi_p}\wedge \alpha_{i_1\dots
i_p}^{\bsi_1\dots\bsi_p},\qquad \alpha_{i_1\dots
  i_p}^{\bsi_1\dots\bsi_p}\in\pi_{r+1,r}^*\left(\Lambda^q_r\right),
\end{equation}
where $\abs{\bsi_l}\leq r$ for $l=1,\dots,p$. Note that:
\begin{enumerate}
\item derivatives of
order $r+1$ appear in the above expression in the forms
$\omega^{i_l}_{\bsi_l}$ with $\abs{\bsi_l}=r$. It is possible to
obtain an expression containing just $r$-th order derivatives by using
contact forms of the type $d\omega^{i_l}_{\bsi_l}$ with
$\abs{\bsi_l}=r-1$; see \cite{Kru90};
\item in the case $q=0$ it can be proved (for $p=1$, see~\cite{Kru90}; in the
  general case the argument is similar, see~\cite{Vit01}) that
  $\abs{\bsi_l}\leq r-1$ for $l=1,\dots,p$ in the above coordinate expression.
\end{enumerate}

From the above consideration it follows that the horizontalisation
allows us to discard contact components from a form. And contact
components produce no contribution to action-like functionals
(section~\ref{sec:var_princ}). Moreover, we shall see that the
first term of the $\cC$-spectral sequence is made by quotients of
$p$-contact forms by $(p+1)$-contact forms (section
\ref{sec:C-ss}), so that it is important to be able to discard
$(p+1)$-contact factors from a $p$-contact form. Hence, for future
purposes, we introduce the following \emph{partial
horizontalisation} map
\begin{equation}\label{eq:5}
h^{p,q}\colon \Lambda ^{p+q}_r \to \Lambda ^{p}_{r+1,r}\otimes
\overline{\Lambda }^{q}_r,\quad
\alpha \mto
(\wedge^p\id\otimes\wedge^q {D^{r+1}}^*)\circ (\pi _{r+1,r}^*\alpha).
\end{equation}
The action of $h^{p,q}$ on decomposable forms is
\begin{multline*}
h^{p,q}(\alpha_{1}\wedge \ldots \wedge \alpha _{p+q})=
\\
\frac{1}{p!\,q!}\sum_{\sigma \in S_{p+q}}\abs{\sigma}
\pi_{r+1,r}^*(\alpha _{\sigma(1)}\wedge \ldots \wedge \alpha
_{\sigma (p)}) \otimes h^{0,q}(\alpha_{\sigma (p+1)}\wedge \ldots
\wedge \alpha _{\sigma (p+q)})\notag ,
\end{multline*}
where $S_{p+q}$ is the set of permutations of $p+q$ elements.
\section{Spectral sequence}
\label{sec:C-ss}

In this section we present a new finite order approach to $\cC$-spectral
sequence on the jets of submanifolds of order $r$.  This approach has been
directly inspired by Vinogradov's $\cC$-spectral sequence on infinite order
jets~\cite{Vin77,Vin78,Vin84}.  A first analysis of the finite order
$\cC$-spectral sequence already appeared in~\cite{MaVi02}. We complete that
exposition with the new coordinate expressions~\eqref{eq:6},
\eqref{pro:coordex_hd} and the relations between our differentials and the
standard horizontal differential $d_H$. We also point out in
remark~\ref{re:hor_coh} that our methods yield a generalisation of results
previously obtained in~\cite{Kru90} and~\cite{Gra00}.

\medskip

The module $\Lambda^k_r$ is filtered by the submodules
$\cC^p\Lambda^{k}_r$; namely, we have the obvious \emph{finite} filtration
\begin{displaymath}
\Lambda^k_r \byd \cC^0\Lambda^{k}_r\supset \cC^1\Lambda^{k}_r \supset
\dots \supset
\cC^p\Lambda^{k}_r \supset \dots \supset \cC^I\Lambda^{k}_r \supset
\cC^{I+1}\Lambda^{k}_r = \{0\},
\end{displaymath}
where $I$ is the codimension of the Cartan distribution
(see~\cite{Many99}). The filtration is stable with respect to the
differential of forms, \emph{i.e.}, $d(\cC^p\Lambda^{k}_r)\subset
\cC^p\Lambda^{k+1}_r$. We say that the above graded filtration is the
\emph{$\cC$-filtration} on the jet space of order $r$.

The $\cC$-filtration gives rise to a spectral sequence $(E^{p,q}_N,
e_N)_{N,p,q \in \mathbb{N}}$ in the usual way.  We recall that a spectral
sequence is a sequence of differential Abelian groups where each term is the
cohomology of the previous one, with the exception of $E^{*,*}_0$, which, in
our case, is the set of quotients between consecutive terms of the
$\cC$-filtration (see, \emph{e.g.},~\cite{KrVe98,Vit99}).

\begin{definition}
  We call the above spectral sequence the \emph{Vinogradov's $\cC$-spectral
    sequence of (finite) order $r$} on $E$.
\end{definition}

Next goal is to recall the description of all terms in the finite
order $\cC$-spectral sequence made in~\cite{MaVi02}.

\medskip

We recall that $E_{0}^{p,q} \equiv \cC^p\Lambda^{p+q}_{r} \big /
\cC^{p+1}\Lambda^{p+q}_{r}$.  Generalising the equalities~\eqref{cC_and_ker_h}
to partial horizontalisation we obtain the following result.
\begin{lemma}
\begin{displaymath}
\cC^{p+1}\Lambda _{r}^{p+q} =\ker\,h^{p,q}\quad\text{if}\quad
q\leq n,\qquad \cC^{p+1}\Lambda _{r}^{p+q} =\Lambda_{r}^{p+q}\quad
\text{if}\quad q>n.
\end{displaymath}
\end{lemma}
As a consequence, we are able to express any equivalence class of
$E_0$ with a distinguished form. To proceed with our
investigation, we need to describe the target space of partial
horizontalisation of contact forms. Taking into account
\eqref{eq:coord_contact}, we introduce the space
$\cC^p\Lambda^p_{r,r+1}\subset\cC^p\Lambda^p_{r+1}$ of contact
forms with coefficients in $\cF_r$. This space can be
characterised in an intrinsic way as follows. A form $\gamma$ is
in $\cC^p\Lambda^p_{r,r+1}$ if and only if
$\gamma=\pi_{r+1,r}^*(i_{X_1}\cdots i_{X_q}\gamma')$, where
$\gamma'\in\cC^p\Lambda^{p+q}_r$ and $X_1,\ldots, X_q\colon
J^1(E,n)\to H^{1,0}$. Of course, the action of $\gamma'$ on the
vector fields is obtained through the isomorphism \eqref{eq:1}.

\begin{proposition}[Computation of $E_0$,~\cite{MaVi02}]\label{iso2}
  Let $q>0$. Then, the restriction of $h^{p,q}$ to $\cC^{p}\Lambda_{r}^{p+q}$
  yields the isomorphism
\begin{displaymath}
E_{0}^{p,q}=\cC^{p}\Lambda
_{r}^{p+q}\big/\cC^{p+1}\Lambda_{r}^{p+q} \to
\cC^{p}\Lambda_{r,r+1}^{p}\otimes \hL_{r}^{q},\,\quad \ [\alpha]
\mto h^{p,q}(\alpha).
\end{displaymath}
\end{proposition}
Obviously, if $q=0$ then $E_0^{0,0}=\cF_r$, and $E_{0}^{p,0}=\cC^{p}\Lambda
_{r}^{p}$ for $p>0$. We set $\hd\byd e_0$; we have
$\hd( h^{p,q}(\alpha ))=h^{p,q+1}(d\alpha)$.
Hence, the bigraded complex $(E_{0},e_{0})$ is isomorphic to the sequence of
complexes $(\cC^{p}\Lambda_{r,r+1}^{p}\otimes
\hL_{r}^{*},\hd)_{p\in\mathbb{N}}$.  The complexes have finite length $n$, and
the sequence is trivial for $p>I$.  For $p=0$ this is just the \emph{horizontal
  de Rham complex of order $r$} (see~\cite{Ver98} for the infinite order
version).

We observe that the coordinate expression of $\bar{\alpha}\in
\cC^{p}\Lambda_{r,r+1}^{p}\otimes \hL_{r}^{q}$ is
\begin{align}
\label{eq:6}
  &\bar{\alpha}= \bar{\alpha}^{\bsi_1\cdots\bsi_p}_{i_1\cdots i_p\
    \lambda_1\cdots\lambda_q}
  \omega^{i_1}_{\bsi_1}\wedge\cdots\wedge\omega^{i_p}_{\bsi_p}
  \otimes\odx^{\lambda_1}\wedge\cdots\wedge\odx^{\lambda_q},
  \\
  &\bar{\alpha}^{\bsi_1\cdots\bsi_p}_{i_1\cdots i_p\;\lambda_1\cdots\lambda_q}=
  \alpha^{\bsi_1\cdots\bsi_p\ \bta_1\cdots\bta_l} _{i_1\cdots i_p\ j_1\cdots
    j_l\;\lambda_{l+1}\cdots\lambda_q} u^{j_1}_{\bta_1,\lambda_1}\cdots
  u^{j_l}_{\bta_l,\lambda_l},\notag
\end{align}
with $0\leq\abs{\bsi_k}\leq r$, $\abs{\bta_h}=r$, $0\leq l\leq q$ and
$\alpha^{\bsi_1\cdots\bsi_p\ \bta_1\cdots\bta_l}
  _{i_1\cdots i_p\ j_1\cdots
    j_l\;\lambda_{l+1}\cdots\lambda_q}\in\cF_r$.

\begin{proposition}\label{pro:coordex_hd}
  Let $\bar\alpha\in \hL_{r}^{q}$.
  Then, we have the coordinate expression
\begin{equation}
\label{eq:8a}
  \hd\bar{\alpha}= D_\lambda(\alpha^{\bta_1\cdots\bta_l}
  _{j_1\cdots j_l\ \lambda_{l+1}\cdots\lambda_q})
  u^{j_1}_{\bta_1,\lambda_1}\cdots u^{j_l}_{\bta_l,\lambda_l}
  \odx^{\lambda}\wedge\odx{}^{\lambda_1}
  \wedge\cdots\wedge\odx{}^{\lambda_{q}}.
\end{equation}
\end{proposition}
\begin{proof}
  A form $\alpha\in\Lambda^q_r$ such that
  $h^{0,q}(\alpha)=\bar{\alpha}$ has the coordinate expression
\[
\alpha=\alpha^{\bta_1\cdots\bta_l}
  _{j_1\cdots j_l\;\lambda_{l+1}\cdots\lambda_q}
du^{j_1}_{\bta_1}\wedge\cdots\wedge du^{j_l}_{\bta_l}\wedge
dx^{\lambda_{l+1}}\wedge\cdots\wedge dx^{\lambda_q}.
\]
From the expression of $d\alpha$ we easily get the result.
\end{proof}

In the case that $\alpha$ belongs to
$\cC^p\Lambda^p_{r,r+1}\otimes\hL^q_r$ for $p\geq 1$, one has also
to differentiate contact forms (see \eqref{eq:6}). We have
$d\omega^i_{\bsi}=-du^i_{\bsi ,\lambda}\wedge dx^\lambda$, hence
the differentiation when $\abs{\bsi}=r$ yields a form on an higher
order jet. So, to produce a coordinate expression for
$\hd\bar\alpha$ we need a form $\alpha$ such that
$h^{p,q}(\alpha)=\bar\alpha$ and whose expression use
$d\omega^i_{\brh}$ with $\abs{\brh}=r-1$, instead of
$\omega^i_{\bsi}$ with $\abs{\bsi}=r$. Such an expression can be
found in~\cite{Kru90} (see also~\eqref{eq:coord_contact}).

\begin{remark}\label{rem:hd-dH}
  It is important to stress that the differential $\hd$ of the previous
  complexes does not coincide with the \emph{horizontal differential} $d_H$
  (denoted by $\widehat{d}$ in~\cite{Many99}; see also~\cite{Sau89}) used in
  the infinite order formalism.  The action of $d_H$ is characterised by the
  coordinate expressions $d_H f=D_\lambda f\, dx{}^\lambda$, for $f\in\cF$,
  $d_H\, dx{}^\lambda=0$, $d_H\omega^i_{\bsi}=-\omega^i_{\bsi,\lambda}\wedge
  dx^{\lambda}$.  Hence, $\hd$ does not change the order of jet space while
  $d_H$ raises the jet order by one.

  Moreover, from the above proposition it turns out that
  $\hd\bar\alpha=d_H\bar\alpha$ if $\bar\alpha$ does not contain derivatives of
  order $r+1$, \emph{i.e.}, if $\bar\alpha\in \cC^{p}\Lambda_{r-1,r}^{p}\otimes
  \hL_{r}^{q}\subset \cC^{p}\Lambda_{r,r+1}^{p}\otimes \hL_{r}^{q}$. In other
  words, $\hd$ and $d_H$ just differ on highest order derivatives. Hence, it is
  clear that these two operators coincide in the direct limit, \emph{i.e.}, on
  infinite order jets.
\end{remark}

We recall that $E_{1} = H(E_{0})$, where the cohomology is taken with respect
to $\hd$. We start by determining the term $E_1^{p,n}$.

\begin{proposition}[\cite{MaVi02}]\label{pro:bicomplex}
  We have the diagram
\newdiagramgrid{finitebicomplex}%
{.8,.8,1,1.6,1.3,.7,.7,1.3,1,1,}{1,.7,1,.7,1,.7,1}
\begin{equation}\label{eq:7}
\begin{diagram}[grid=finitebicomplex]
  & & 0 & & 0 & & & & 0 & &
  \\
  & & \uTo & & \uTo & & & & \uTo & &
  \\
  0 & \rTo & \hL^n_r \big/ \hd(\hL^{n-1}_r) & \rTo^{e_1} & E_0^{1,n}
  /\hd(E_{0}^{1,n-1}) & \rTo^{e_1} & \dots & \rTo^{e_1} & E_0^{p,n}
  /\hd(E_{0}^{p,n-1}) & \rTo^{e_1} & \dots
  \\
  & & \uTo & \ruTo(2,2)_{\widetilde{e_1}} & \uTo & & & & \uTo & &
  \\
  & & \hL^n_r & & E_0^{1,n} & & \cdots & & E_0^{p,n} & & \cdots
  \\
  & & \uTo_{\hd} & & \uTo_{\hd} & & & & \uTo_{\hd} & &
  \\
  & & \cdots & & \cdots & & \cdots & & \cdots & & \cdots
\end{diagram}
\end{equation}
where $E_1^{*,n}$ is the top row, and
$e_1([h^{p,n}(\alpha)])=[h^{p+1,n}(d\alpha )]$.
\end{proposition}
For a proof, see~\cite{MaVi02}. In the above diagram, the maps without a label
are either trivial of quotient projections. The map $\widetilde{e_1}$ is just
the composition of $e_1$ with the quotient projection $\hL^n_r\to\hL^n_r\big
/\hd(\hL^{n-1}_r)$.

\begin{remark}
In order to match a similar convention in the infinite order
formulation, in \eqref{eq:7} we can replace $\hd$ with
$(-1)^p\hd$.
\end{remark}

\medskip

In order to determine $E_1^{p,q}$ with $q<n$ we need some preliminary results.
First of all, the sequence
\newdiagramgrid{contact}
{.8,.8,.8,.8,.8,.8,.8,.8,.8,.8,.8} {1}
\begin{equation}\label{eq:contact_seq}
\begin{diagram}[grid=contact]
  0 & \rTo & \cC^p\Lambda_r^{p} & \rTo^d & \cC^p\Lambda_r^{p+1} & \rTo^d &
  \dots & \rTo^d & \cC^p\Lambda_r^{p+n-1} & \rTo^d & \dots
\end{diagram}
\end{equation}
is exact up to the term $\cC^p\Lambda_r^{p+n-1}$~\cite{MaVi02}. The
$\cC$-spectral sequence converges to the de Rham cohomology of $J^r(E,n)$,
because it is a first quadrant spectral sequence~\cite{MaVi02}. Moreover, the
cohomology of $J^r(E,n)$ is equal to the cohomology of $J^1(E,n)$ due to the
topological triviality of the fibres.  Summing up these facts, we have
\begin{enumerate}
\item $E_1^{0,q}=H^q(J^1(E,n))$, for $q\neq n$;
\item $E_2^{p,n}=H^p(J^1(E,n))$, for $p \geq 1$;
\item $E_1^{p,q}=0$ for $q\neq n$ and $p\neq 0$.
\end{enumerate}
This completes the description of all terms of the finite order $\cC$-spectral
sequence. The interested reader can find details and proofs in~\cite{MaVi02}.
\begin{remark}\label{re:hor_coh}
  A partial result in the direction of the computation of the finite order
  $\cC$-spectral sequence on jets of submanifolds already appeared
  in~\cite{Gra00}. In that paper the term $E^{0,q}_1$ was computed with
  interesting intrinsic techniques (Koszul complexes), by analogy with the
  infinite order case. This term is known as \emph{horizontal
    cohomology}~\cite{Ver98}.
  
  However, we compute \emph{all} terms of the finite order $\cC$-spectral
  sequence. Our method is a direct generalisation of the argument that Krupka
  used to compute $E^{0,q}_1$ in his setting on jets of fibrings~\cite{Kru90}.
  This is a direct consequence of the exactness of the
  sequence~\eqref{eq:contact_seq} for $p=1$, which is proved in~\cite{Kru90}
  through a `vertical' Poincar\'e lemma and elementary facts from sheaf theory,
  from which the results in~\cite{Gra00} follow.
\end{remark}
\section{Green-Vinogradov formula and the finite order variational sequence}
\label{sec:green-vinogr-form-1}

The diagram~\eqref{eq:7} contains a further complex which is of fundamental
importance for the calculus of variations: the variational sequence.

\begin{definition}\label{def:varseq}
  The complex
  \newdiagramgrid{varcomplex} {.7,.7,1,.7,1,1.3,2,1,2,.5,.7}{1}
\begin{displaymath}
\begin{diagram}[grid=varcomplex]
  \dots & \rTo^{\hd} & \hL^{n-1}_r & \rTo^{\hd} & \hL^{n}_r &
  \rTo^{\Tilde{e}_1} & E_0^{1,n} /\hd(E_{0}^{1,n-1}) & \rTo^{e_1}
  & E_0^{2,n} /\hd(E_{0}^{2,n-1}) & \rTo^{e_1} & \dots,
\end{diagram}
\end{displaymath}
is said to be the \emph{variational sequence of order $r$}.

Elements of $E_0^{p,n} /\hd(E_{0}^{p,n-1})$ are said to be
\emph{variational $p$-forms of order $r$}.
\end{definition}
From the calculations of the terms of the $\cC$-spectral sequence it follows
that the cohomology of the above complex is isomorphic to the de Rham
cohomology of $J^1(E,n)$.

In the last section we shall analyse the relationship between the
variational sequence and the calculus of variations; for the
moment, we concentrate on another problem. Namely, variational
forms are equivalence classes; we are going to prove that each
equivalence class can be represented by a distinguished form.  To
this purpose, the most important tool is the Green--Vinogradov
formula, which is the geometric analogue of the integration by
parts.
\subsection{Green--Vinogradov formula}
\label{ssec:green-vinogr-form}

In this subsection we first exhibit a natural isomorphism between
the module of contact forms and a suitable space of differential
operators. This allows us to `import' the theory of
adjoint operators and the Green--Vinogradov formula in our setting
(see~\cite{Many99,KrVe98} for introductory expositions).

Here, we also show that this procedure does not have a finite
order analogue. So, we are forced to use the infinite order theory
to provide distinguished representatives for equivalence classes
in the variational sequence. This is well-known in the infinite
order case, but its application to the finite-order case is new in
the case of jets of submanifolds, up to partial results contained
in~\cite{MaVi02}.

We observe that our representation solves a problem which was left open
in~\cite{Gra00} (remark~\ref{re:griff_conj}).

\medskip

Let $P$, $Q$ be projective modules over an $\mathbb{R}$-algebra
$A$. We recall (see~\cite{AVL91} for instance) that a \emph{linear
differential operator} of order $k$ is an $\mathbb{R}$-linear map
$\Delta \colon P\to Q$ such that
\begin{equation*}
[\delta _{a_{0}},[\dots ,[\delta _{a_{k}},\Delta ]\dots]]=0
\end{equation*}
for all $a_{0},\dots ,a_{k}\in A$. Here, square brackets stand for
commutators and $\delta _{a_{i}}$ is the multiplication morphism
by $a_{i}$. Of course, linear differential operators of order zero
are morphisms of modules. This definition can be generalised to
differential operators in several arguments.

We are interested in differential operators between
$\cF_r$-modules whose expressions contain total derivatives
instead of partial ones.  Unfortunately, such operators raise the
jet order by one, so we are lead to extend the above definition in
an obvious way to the case in which $P$ is an $\cF_{r}$-module and
$Q$ is an $\cF_{s}$-module, with $r\leq s$ (taking into account
the inclusion $\cF_{r}\subset \cF_{s}$).

Now, let $k\leq s-r$; we say that a differential operator $\Delta
\colon P\to Q$ of order $k$ is
\emph{$\mathcal{C}$-dif\-fer\-en\-tial} if it can be restricted to
manifolds of the form $L^{(r)}$ and $L^{(s)}$. In local
coordinates, $\mathcal{C}$-dif\-fer\-en\-tial operators have the
matrix form $(a_{ij}^{\bsi}D_{\bsi})$, where $a_{ij}^{\bsi}\in
\cF_{s}$, $D_{\bsi}=D_{\sigma _{1}}\circ \dots \circ D_{\sigma
_{h}}$, and $\abs{\bsi}=h\leq k$.  The space of such operators is
denoted by $\CDiff_k(P,Q)$.  We shall deal with spaces of
antisymmetric $\mathcal{C}$-differential operators of $l$
arguments in $P$, which we denote by
$\CDiff^{\alt}_{(l)\,k}(P,Q)$.

\medskip

Next, we consider sections of pseudo-vertical bundles.  Let $X$ be a vector
field on $J^r(E,n)$. Its \emph{vertical part} $X_V$ is the
section of $V^{r+1,r}$ defined by $X_V\byd\omega^{r+1}(X\circ\pi_{r+1,r})$.  If
$X=X^\lambda \pd{}{x^\lambda} + X^i_{\bsi}\pd{}{u^i_{\bsi}}$, then
$X_V=(X^i_{\bsi}-u^i_{\bsi,\lambda}X^\lambda)[\pd{}{u^i_{\bsi}}]$.

For $r\geq 1$ let us denote the $\cF_{r}$-module of bundle
morphisms $\varphi \colon J^{r}(E,n)\to V^{1,0}$ over
$\id_{J^{1}(E,n)}$ by $\varkappa _{r}$.  Of course $\varkappa
_{1}$ is just the set of the sections of $V^{1,0}$. We also define
$\varkappa_0\subset\varkappa_1$ as the subset of vertical parts of
vector fields on $E$.  Any $\varphi\in\varkappa_{r}$ can be
uniquely prolonged to an \emph{evolutionary vector field}
$\Evo_{\varphi }\colon J^{r+s}(E,n)\to V^{s+1,s}$
(see~\cite{Many99}). In coordinates, if $\varphi
=\varphi^{i}\,[\pd{}{u^{i}}]$, then
$\Evo_{\varphi}=D_{\bsi}\varphi^{i}\,[\pd{}{u_{\bsi}^{i}}]$, with
$\abs{\bsi}\leq s$.
\begin{proposition}[\cite{MaVi02,Vit01}]\label{iso1}
  We have the natural isomorphism
\begin{displaymath}
  \cC^{p}\Lambda _{r,r+1}^{p}\otimes \bar{\Lambda}_{r}^{q} \to
  \CDiffalt{p}{r}(\varkappa _{0},\bar{\Lambda}_{r}^{q}),\, \bar{\alpha} \mto
  \nabla _{\bar{\alpha}}(\varphi_{1},\ldots ,\varphi_{p})= \frac{1}{p!}\;
  i_{\Evo_{\varphi_1}}(\cdots i_{\Evo_{\varphi_p}}(\bar{\alpha})\cdots ).
\end{displaymath}
\end{proposition}

The above proposition can be proved by analogy with the infinite order case
(see~\cite{Many99}). The numerical factor is put in
order to have simpler coordinate expressions: from \eqref{eq:6} we have
\begin{equation}\label{eq:8}
  \nabla_{\bar{\alpha}}(\varphi_{1},\ldots ,\varphi_{p})=
   \alpha^{\bsi_1\cdots\bsi_p}_{i_1\cdots
  i_p\;\;\lambda_1\cdots\lambda_q} D_{\bsi_1}\varphi^{i_1}_1\cdots
  D_{\bsi_p}\varphi^{i_p}_p\,
\odx^{\lambda_1}\wedge\ldots\wedge\odx^{\lambda_q}.
\end{equation}

\medskip

Let $\cF$ be the direct limit of the chain
$\cdots\subset\cF_r\subset\cF_{r+1}\subset\cdots$ of injections given by
pull-back. In other words, $\cF$ is the space of all functions on any finite
order jet. In what follows, we shall drop the order index $r$ to indicate
spaces obtained in similar ways.  We say that $P$ is a \emph{horizontal module}
if $P$ is an $\cF$-module obtained as direct limit of an ascending chain of
injections of $\cF_r$-modules.  As an example, the injections $\varkappa
_{i}\hookrightarrow \varkappa_{i+1}$ yield a sequence whose direct limit,
denoted by $\varkappa$, is a horizontal module. See~\cite{Many99} for
connections between $\varkappa$ and symmetries of differential equations.

Let $p>0$; for any horizontal module $P$ we have the complexes
$(\CDiff_{(p)}(P,\hL^*),w)$, where $w(\nabla)=d_H\circ\nabla$. In~\cite{Many99}
it is proved that the only non-vanishing cohomology group of such complexes is
the $n$-th, which is equal to $\CDiff_{(p-1)}(P,\widehat{P})$, where
$\widehat{P}\byd\Hom_\cF(P,\hL^n)$.  Any $\cC$-differential operator
$\Delta\colon P\to Q$ induces a map of the corresponding $w$-complexes, hence a
cohomology map $\Delta^*\colon \widehat{Q}\to\widehat{P}$. It fulfils the
\emph{Green-Vinogradov} formula~\cite{Vin77,Vin78}
\begin{equation}
  \label{eq:adjoint}
  \widehat{q}(\Delta(p))-(\Delta^*(\widehat{q}))(p) =
  d_H\omega_{p,\widehat{q}}(\Delta).
\end{equation}
In coordinates, if $\Delta=\Delta_{ij}^{\bsi}D_{\bsi}$, then
$\Delta^*=(-1)^{\abs{\bsi}}D_{\bsi}\circ\Delta_{ji}^{\bsi}$. One can easily
realize the meaning of the above formula by observing that the left-hand side
contains no zero order term with respect to $p$, hence it must be the total
derivative of a certain $\omega_{p,\widehat{q}}$.  The problem of representing
variational forms in the infinite order formalism is solved by taking the
skew-symmetric part of the cohomology of the complexes
$(\CDiff_{(p)}(P,\hL^*),w)$ with respect to permutations of the arguments
(\cite{Vin77,Vin78,Vin84}; see also~\cite[p. 192]{Many99}).  Namely, let
$K_p(\varkappa)\subset\CDiff^{\alt}_{(p-1)} (\varkappa,\widehat{\varkappa})$ be
the subspace of operators $\nabla$ fulfilling
$\nabla(\varphi_1,\ldots,\varphi_{p-2})^*=
-\nabla(\varphi_1,\ldots,\varphi_{p-2})$ for all $\varphi_i\in\varkappa$,
$i=1,\ldots,p$. Then, we have the isomorphism
\begin{equation}\label{eq:repr}
  I_p\colon \CDiffalt{p}{}(\varkappa, \hL^{n})\big /
     \hd(\CDiffalt{p}{}(\varkappa,\hL^{n-1})) \to
     K_p(\varkappa)\,,\quad
[\Delta]\mto \overline{\Delta}.
\end{equation}
where $\overline\Delta(\varphi_{1},\ldots,\varphi_{p-1})(\varphi_p)\byd
\Delta(\varphi_{1},\ldots,\varphi_{p-1})^*(1)(\varphi_p)$
for all $\varphi_i\in\varkappa$, $i=1,\ldots,p$.

\medskip

Now we devote ourselves to the problem of the representation of the finite
order variational sequence. We would like to find representatives of the same
order of all objects in the given class. But the following attempts, carried
out by analogy with the infinite order case, failed.
\begin{enumerate}
\item We would like to reproduce the above scheme in the finite order case. We
  could generalise finite order $\hd$-complexes of antisymmetric operators to
  finite order complexes with arguments in a module $P$ over $\cF_r$ of
  operators of any kind.  Unfortunately, any $\cC$-differential operator
  $\Delta\colon P\to Q$ of order $k>0$ is a \emph{graded} map of such finite
  order complexes, in the sense that it raises the order by $k$. Hence, there
  is no hope to find a version of `adjoint operator' which is order-preserving.
\item Let $p=1$ for simplicity. Given $\Delta\in \CDiffalt{1}{r}(\varkappa_0,
  \hL^{n}_r)$ we can look for operators
  $\bar{\omega}_\Delta\in\CDiffalt{1}{r}(\varkappa_0,\hL^{n-1}_r)$ and
  $\bar\Delta\in\Hom_{\cF_r}(\varkappa_0,\hL^n_r)$ such that
  $\Delta+\hd(\bar{\omega}_\Delta)=\bar\Delta$. Unfortunately, the previous
  equation always admits a unique solution if we let $\bar\Delta$ and
  $\bar{\omega}_\Delta$ to be defined on an higher order jet. A proof of this
  statement can be easily achieved by the same argument used by
  Kol\'a\v r~\cite{Kol83} and the coordinate expression of $\hd$
  (proposition~\ref{pro:coordex_hd}).
\end{enumerate}

The problem of representation is thus solved by embedding the finite order
variational sequence into the infinite order variational sequence and using the
standard results of the infinite order theory restricted to the image of the
embedding.  To this aim, we observe that pull-back allows us to take $\hd=d_H$
(remark~\ref{rem:hd-dH}). Moreover, pull-back is a map of complexes between
$\cC$-variational sequences of different order~\cite{Vit99}, hence the direct
limit of the $r$-th order $\cC$-variational sequence is just the standard
infinite order $\cC$-variational sequence.
\begin{theorem}\label{th:iso}
  The representation of each space of finite order variational sequence (in the
  case $p>1$) is the image space $K_{p,r}(\varkappa_0)\byd
  \im(I_p\circ\chi_{p,r})$, where
\begin{multline*}
  \chi_{p,r}\colon\CDiffalt{p}{r}(\varkappa_0, \hL^{n}_r)\big /
  \hd(\CDiffalt{p}{r}(\varkappa_0,\hL^{n-1}_r)) \hto
  \\
  \CDiffalt{p}{}(\varkappa, \hL^{n})\big /
  \hd(\CDiffalt{p}{}(\varkappa,\hL^{n-1}))
\end{multline*}
is the inclusion induced by pull-back $\pi_{r+1,r}^*$ into the direct limit.
\end{theorem}
\begin{remark}\label{re:griff_conj}
  In the paper~\cite{Gra00} there is a \emph{conjecture/question} by Griffiths
  about the existence of natural representatives for $E^{1,n}_1$ in the setting
  of finite order jets of submanifolds. The above theorem answers
  constructively to this question. We also provide the coordinate expression of
  such representatives in the next subsection.
\end{remark}

\begin{remark}
  We could consider the `complementary' problem to the representative's one.
  More precisely, given $\alpha\in\CDiffalt{p}{r}(\varkappa_0, \hL^{n}_r)$ we
  can look for an operator $\omega$ fulfilling $\alpha = \nabla_\alpha + \hd
  \omega$.  Such a section always exists due to the vanishing of the cohomology
  of $\hd$ on the space where $\omega$ lives.

  This problem, for $p=1$, amounts to search a Poincar\'e--Cartan form
  (see~\cite{Kol83} for a discussion and its solution). In the case of jets of
  fibrings, a form of that kind can be determined through a linear symmetric
  connection on the base manifold. See~\cite{AB01} for the general situation
  $p>1$.  We expect that the above results hold also in the more general
  framework of jets of submanifolds. See~\cite{Gri97,Gri98} for a first
  discussion of the problem.
\end{remark}
\begin{remark}
  In the case of jets of fibrings the problem of the representation has been
  faced within Krupka's formulation of variational sequences too.  We
  mention~\cite{Vit98}, with results up to $p=n+2$, and~\cite{KrMu03}.  But
  Krupka's formulation produces the same variational sequence as the finite
  order $\cC$-spectral sequence~\cite{Vit99}, where the representation problem
  has been considered in~\cite{Vit01} and solved with similar techniques as
  here.
\end{remark}
\subsection{Coordinate expressions}
\label{ssec:coord-expr}

Here we shall provide coordinate expressions for the representative of
variational forms that we found in theorem~\ref{th:iso}. We stress that such
expressions have never been written before for $p\geq 3$ in the case of jets of
submanifolds. As a by-product, this provides the well-known interpretation of
the variational sequence in terms of the calculus of variations
(section~\ref{sec:var_princ}). We also obtain the coordinate expressions of the
differentials of the variational sequence.

\medskip

An operator $\Delta\in\CDiffalt{p}{r}(\varkappa_0,\hL^{n}_r)$ has the
expression (see also \eqref{eq:6})
\begin{equation}
  \Delta(\varphi_1,\ldots,\varphi_p) =
   \Delta_{i_1\dots i_{p-1}\,j}^{\bsi_1\dots\bsi_{p-1}\bta}
  D_{\bsi_1}\varphi^{i_1}_1\cdots D_{\bsi_{p-1}}\varphi^{i_{p-1}}_{p-1}
  D_{\bta}\varphi^j_p\,\Vol_n,
\end{equation}
where $0\leq\abs{\bsi_a}\leq r$, $0\leq\abs{\bta}\leq r$,
$a=1,\ldots,p-1$, $\Delta_{i_1\dots
i_{p-1}\,j}^{\bsi_1\dots\bsi_{p-1}\bta}$ is a polynomial with
respect to $(r+1)$-st derivatives of distinguished type (it
contains hyperjacobians: see~\eqref{eq:4} and related comments,
remark~\ref{re:hyperjac} and~\cite{PaVi00}) with coefficients in
$\cF_r$ and $\Vol_n \byd n!\,\odx^1\wedge\dots\wedge\odx^n$ is a
local volume form on any submanifold $L\subset E$ which is
concordant with the given chart. We shall use the following
Leibnitz formula for total derivatives (see~\cite{Sau89} for the
case $k=2$)
\[
D_{\bsi}(f_1\cdots f_k)=
\sum_{\abs{\bsi_{1}}+\abs{\bsi_{2}}+\cdots+
\abs{\bsi_{k}}=\abs{\bsi}}
\frac{\bsi!}{\bsi_{1}!\bsi_{2}!\cdots\bsi_{k}!}
D_{\bsi_1}f_1\cdots D_{\bsi_k}f_k.
\]
We also denote by $(\bta_1,\bta_2)$ the multiindex which is the union of
$\bta_1$ and $\bta_2$.

\medskip

In the case $p=1$ we have $I_1([\Delta])(\varphi_1) =
\overline{\Delta}(\varphi_1)= (-1)^{\abs{\bta}}D_{\bta}\Delta_i^{\bta}
\varphi^i_1\,\Vol_n$.

\medskip

In the case $p=2$ we have
\begin{align}\label{eq:I2}
  \overline{\Delta}(\varphi_1)(\varphi_2) &=
  (-1)^{|\bta|}D_{\bta}\left(\Delta_{i_1\
  j}^{\bsi_1\bta}
  D_{\bsi_1}\varphi^{i_1}_1\right)\,\varphi^j_2 \, \Vol_n
  \\
  &= \sum_{0\leq |\bnu|+|\bta_1|\leq
  r}(-1)^{\abs{(\bnu,\bta_1)}}\frac{(\bnu,\bta_1)!}{\bnu!\bta_1!}
  D_{\bnu}\Delta_{i_1}^{\bsi_1}{}_j^{(\bnu,\bta_1)}
  D_{(\bta_1,\bsi_1)}\varphi^{i_1}_1\varphi^{j}_2 \, \Vol_n \notag
  \\
  &= \left(\sum_{\substack{(\bsi_1,\bta_1)=\brh_1\\ 0\leq |\brh_1|\leq 2r}}
    (-1)^{\abs{(\bnu,\bta_1)}}\frac{(\bnu,\bta_1)!}{\bnu!\bta_1!}
    D_{\bnu}\Delta_{i_1\ j}^{\bsi_1(\bnu,\bta_1)}\right)
  D_{\brh_1}\varphi^{i_1}_1 \varphi^{j}_2 \, \Vol_n, \notag
\end{align}
where the last passage follows after renaming multiindexes and rearranging
terms.

\medskip

In the case $p\geq 3$ we proceed by analogy with the case $p=2$ with only the
change of the length of the product on which to apply Leibnitz rule. We have
\begin{alignat}{2}\label{eq:Ip}
  \overline{\Delta} (\varphi_1 &,\ldots,\varphi_{p-1})(\varphi_p) =
\\
  = &
  (-1)^{|\bta|}D_{\bta}(\Delta_{i_1\cdots i_{p-1}\ j}
   ^{\bsi_1\cdots\bsi_{p-1}\bta}
  D_{\bsi_1}\varphi^{i}_1\cdots D_{\bsi_{p-1}}\varphi^{i_{p-1}}_{p-1})
  \varphi^j_p \, \Vol_n \notag
\\
  = &\sum_{0\leq|\bnu|+|\bta_1|+\cdots+|\bta_{p-1}|\leq
   r}(-1)^{\abs{(\bnu,\bta_1,\ldots,\bta_{p-1})}}
  \frac{(\bnu,\bta_1,\ldots,\bta_{p-1})!}{\bnu!\bta_1!\cdots\bta_{p-1}!}
  D_{\bnu}\Delta_{i_1\cdots i_{p-1}\ j}
   ^{\bsi_1\cdots\bsi_{p-1}(\bnu,\bta_1,\ldots,\bta_{p-1})} \times \notag
\\
  & \qquad\times
  D_{(\bta_1,\bsi_1)}\varphi^{i_1}_1\cdots D_{(\bta_{p-1},\bsi_{p-1})}
  \varphi^{i_{p-1}}_{p-1}\varphi^{j}_p \, \Vol_n \notag
\\
  = &\left(\sum_{\substack{(\bsi_i,\bta_i)=\brh_i,\ 0\leq|\brh_i|\leq 2r,
\\[1mm] 1\leq i\leq p-1}}
    (-1)^{\abs{(\bnu,\bta_1,\ldots,\bta_{p-1})}}
  \frac{(\bnu,\bta_1,\ldots,\bta_{p-1})!}{\bnu!\bta_1!\cdots\bta_{p-1}!}
  D_{\bnu}\Delta_{i_1\cdots i_{p-1}\ j}
    ^{\bsi_1\cdots\bsi_{p-1}(\bnu,\bta_1,\ldots,\bta_{p-1})}\right)\times
  \notag
\\
  & \qquad\times D_{\brh_1}\varphi^{i_1}_1\cdots D_{\brh_{p-1}}
    \varphi^{i_{p-1}}_{p-1}\varphi^{j}_p \, \Vol_n. \notag
\end{alignat}

Now, we shall derive from proposition~\ref{pro:bicomplex} the
expression of the differentials $\tilde{e}_1$ and $e_p$, $p>1$
(the expression of $\hd$ has been derived in
proposition~\ref{pro:coordex_hd}).
\begin{remark}
  In general, given a form $\nabla\in K_{p,r}(\varkappa_0)$ it is rather
  difficult to find a form $\alpha\in\Lambda^{p+q}_r$ such that
  $\nabla=I_p([h^{p,n}(\alpha)])$. The commutativity of inclusions between
  variational sequences of different orders implies that, locally,
  $e_p(\nabla_{\bar{\alpha}})=[d\bar{\alpha}]$ (with the symbols of proposition
  \ref{iso1}): this last expression is much easier to be computed. Note that
  $d\bar\alpha$ has to be meant as follows: consider~\eqref{eq:6} as a form on
  a jet space, removing restriction bars, then take the ordinary differential.
  Of course, the resulting expression has an intrinsic meaning.
\end{remark}

If $\lambda\in\hL^n_r$ has the coordinate expression $\lambda=\lambda_0\Vol_n$,
then
\begin{equation}\label{eq:E-L}
\tilde{e}_1(\lambda)(\varphi_1)=
(-1)^{\abs{\bta}}D_{\bta}\left(\pd{}{u^i_{\bta}}\lambda_0\right)
\, \varphi^i_1\Vol_n.
\end{equation}
If $\nabla_{\bar{\alpha}}\in K_{p,r}(\varkappa_0)$, then
\begin{align}
\label{eq:e_p}
e_p(\nabla_{\bar{\alpha}})(\varphi_1,\ldots,\varphi_{p+1}) & =
I_p(d\bar{\alpha})(\varphi_1,\ldots,\varphi_{p+1})
\\
\notag
& = I_p\left(\pd{}{u^{i_1}_{\bsi_1}}
\nabla_{i_2\dots i_{p+1}}^{\bsi_2\dots\bsi_{p+1}}
D_{\bsi_1}\varphi^{i_1}_1\cdots D_{\bsi_{p+1}}\varphi^{i_{p+1}}_{p+1}
\, \Vol_n \right).
\end{align}
The coordinate expression of $e_p$ is obtained by a straightforward
substitution of the coefficients $\pd{}{u^{i_1}_{\bsi_1}}
\Delta_{i_2\dots i_p\,j}^{\bsi_2\dots\bsi_p}$ into \eqref{eq:Ip}.
\section{Variational principles and $\cC$-spectral sequence}
\label{sec:var_princ}

In this section we develop a formalism for the calculus of variations in an
intrinsic geometrical setting, in the general case of $n$ independent
variables. Such a construction already exists for jets of sections (see
\cite{KMS93,Kru01,Kup80,Sau89,Tul77}, for example).

In the case of jets of submanifolds, we found mainly two approaches in
literature: the \emph{parametric} approach (see
\cite{CrSa03a,CrSa03b,Gri97,Gri98} for the general case of $n$ arbitrary; the
approach dates back to Carath\'eodory and earlier for the case $n=1$
\cite{GiHi96}) and Dedecker's approach \cite{Ded77}.  In the parametric
approach the variational principle is formulated on the space $\Imm
J^r_0(\R^n,E)$ under the hypothesis that the Lagrangian commute with the action
of the group of parametrisations (see subsection \ref{ssec:jsubm_top}).  This
leads to extra computations in order to verify at each step the invariance of
objects with respect to changes of parametrisation.

Dedecker tried to reproduce the approach on jets of fibrings, but
he was forced to use families of Lagrangians defined on open
subsets with the property that, on intersecting subsets, the
action be the same.  This considerably complicated his formalism.
Instead we are able to use single objects as Lagrangians because
of the introduction of the pseudo-horizontal bundle and
horizontalisation.

We formulate variational principles on jets of submanifolds in a way as close
as possible to the case of jets of fibrings
\cite{Kru01,Kup80,Sau89,Tul77,Vit98} (which is a close rephrasing in geometric
terms of the standard variational principle). The result is very clean and
provides an interpretation of the variational sequence in terms of calculus of
variations.

\bigskip

\begin{definition}\label{def:genlag}
  A form $[\alpha]=h^{0,n}(\alpha)=\lambda \in\hL_r^n$ is said to be an
\emph{$r$-th order generalised Lagrangian}.
\end{definition}
Indeed, $\lambda$ depends on $(r+1)$-st derivatives in the way specified in
equation~\eqref{eq:4} (\emph{i.e.}, through hyperjacobians;
see~\cite{Olv83,PaVi00}).

\begin{definition}\label{def:action}
The \emph{action} of the Lagrangian $\lambda$ on an $n$-dimensional oriented
submanifold $L\subset E$ with compact closure and regular boundary is the real
number
\begin{equation}\label{eq:11}
\mathcal{A}_L(\lambda)\byd\int_{L}(j_{r}L)^{\ast}\alpha.
\end{equation}
\end{definition}
Due to lemma~\ref{lem:app}, only the horizontal part of a form $\alpha$
contributes to the action, so that the action itself is well-defined.
\begin{remark}
  We are able to introduce a distinguished Lagrangian form; it represents the
  whole class of its Lepage equivalents~\cite{Ded77}. In~\cite{CrSa03a,CrSa03b}
  instead, the authors introduce a distinguished Lepage equivalent called
  Hilbert--Carath\'eodory form. However, it appears that such a form can be
  introduce only if the Lagrangian form is non-vanishing. We do not need such
  an hypothesis.
\end{remark}

\medskip

Now, we formulate the variational problem, \emph{i.e.}, the problem of finding
extremals of the action. Let $L\subset E$ be as in the above definition.  A
vector field $X$ on $E$ vanishing on $\partial L$ is said a
\emph{variation field}.  The submanifold $L$ is \emph{critical} if for each
variation field $X$ with flow $\phi_t$ we have
\begin{equation}
\frac{d}{dt}\Big |_{t=0} \int_{L}(J_{r}\phi_{t}\circ j_{r}L)^{\ast}\alpha =0
\end{equation}
where $J_{r}\phi_{t}:J^{r}(E,n)\longrightarrow J^{r}(E,n)$ is the jet
prolongation of $\phi_{t}$ (see subsection~\ref{ssec:jets}).

We shall show that, indeed, the above condition depends on the vertical part
$X_V$ of $X$, and provide the Euler--Lagrange equations.  First of all, we
observe that ${X_r}_V=\Evo_{X_V}$~\cite{Many99}. We have
\begin{align}
\frac{d}{dt}\Big |_{t=0} \int_{L}(J_{r}\phi_{t}\circ j_{r}L)^{\ast}\alpha
 &=\int_{L}(j_{r}L)^*\cL_{X_{r}}\alpha\label{eq:12}
\\
 &=\int_{L}(j_{r}L)^*i_{X_{r}}d\alpha\label{eq:13}
\\
 &=\int_{L}(j_{r+1}L)^*i_{\Evo_{X_V}}h^{1,n}(d\alpha)\label{eq:14}
\\
 &=\int_{L}(j_{2r+1}L)^*i_{X_V}\tilde{e}_1(\lambda).\label{eq:15}
\end{align}
Here, $\cL_{X_{r}}$ stands for Lie derivative. Equation~\eqref{eq:13} comes
from Stokes'theorem and $\cL_{X_{r}}=i_{X_{r}}d+di_{X_{r}}$.
Equation~\eqref{eq:14} comes from lemma~\ref{lem:app} and the identity
$h(i_{X_{r}}d\alpha)=i_{\Evo_{X_V}}h^{1,n}(d\alpha)$, which is a direct
consequence of the definition of $h^{1,n}$~\eqref{eq:5}.  Finally,
equation~\eqref{eq:15} comes from the identities $i_{\Evo_{X_V}}\hd\beta=\hd
i_{\Evo_{X_V}}\beta$ and $(j_{r+1}L)^*\hd\lambda=d(j_{r+1}L)^*\lambda$, for
$\beta\in E^{1,q}_0$, implying that the value of the integral depends on the
value of the $\hd$-cohomology class of $h^{1,n}(d\alpha)$.

By virtue of the fundamental lemma of calculus of variations,
equation~\eqref{eq:15} vanishes if and only if
\begin{equation*}
(j_{2r+1}L)^*\tilde{e}_1(\lambda)=0,
\end{equation*}
or, that is the same $\tilde{e}_1(\lambda)\circ j_{2r}L=0$.
\begin{remark}
  We obtain intrinsic Euler--Lagrange equation in our scheme, for an arbitrary
  number of independent variables and order of Lagrangian. To do this, we
  followed closely the case of jets of fibrings, hence our construction is
  very natural~\cite{Kru01,Kup80,Sau89,Tul77,Vit98}.

  Obtaining the same in the parametric framework~\cite{CrSa03a,CrSa03b} is much
  more involved. One of the problems is that the authors must prove at each
  step that the objects that they compute have the required invariance with
  respect to the change of parametrisation. Instead, our objects are always
  parametrisation-independent.
\end{remark}

The interpretation of the variational sequence in terms of calculus of
variations is now clear:
\begin{itemize}
\item $\hL^n_r$ is the space of \emph{Lagrangians},
\item $E_1^{1,n}$ is the space of \emph{Euler--Lagrange type forms},
\item $E_1^{2,n}$ is the space of \emph{Helmholtz--Sonin type forms};
\item $\tilde{e}_1$ takes a Lagrangian into its Euler--Lagrange form, the
  vanishing of $\tilde{e}_1$ implies that the Lagrangian is trivial (or
  \emph{null}, see~\cite{And,Olv83});
\item $e_1$ takes an Euler--Lagrange type form into its Helmholtz-Sonin form,
the  vanishing of ${e}_1$ implies that the Euler--Lagrange type form comes from
a Lagrangian.
\end{itemize}
There is no interpretation in terms of known quantities from the calculus of
variations of variational forms for $p\geq 3$ and their differential for $p\geq
2$.

We think that the above variational formalism provides one of the main
motivations for the $\cC$-spectral sequence itself.
\begin{remark}
  We can state what is the form of the most general null Lagrangian $\lambda
  \in\hL_r^n$: locally, it is of the form $\hd p$, where $p\in\hL^{n-1}_r$. Its
  coordinate expression is~\eqref{eq:8a}, with $q=n$.

  Similar statements can be made at any point of the finite order variational
  sequence, and answer very cleanly to questions like the ones considered
in~\cite{CrSa03c} (namely, the structure of null Lagrangians in the parametric
  formalism).

  For example, a \emph{symplectic operator} is an element $B\in E_1^{2,n}$ such
  that $e_1(B)=0$~\cite{Many99} (also known as \emph{Dirac structure}
  in~\cite{Dor93}). The local exactness of the finite order variational
  sequence ensures the local existence of a potential of $B$ in $ E_1^{1,n}$ of
  a definite order and a definite structure of its coefficients.
\end{remark}
\begin{remark}
  Despite the fact that we obtained the above intrinsic theory, we regret that
  we still cannot support with interesting examples, like the ones
  in~\cite{CrSa03a}, the previous considerations, for reasons of time and
  space.  This will be the subject of forthcoming work. Partial results have
  been collected in~\cite{Man03}.
\end{remark}
\section{Conclusions}

This work completes our research on the geometry of jets of submanifolds and
their finite order $\cC$-spectral sequence. But there is one important step
further in finite order theories: the study of the constrained case,
\emph{i.e.}, to repeat the above analysis when the given space is no longer
$J^r(E,n)$ but a submanifold $\mathcal{E}\subset J^r(E,n)$. This will be the
subject of future research.

Moreover, we think that our research constitutes a starting point for a
geometric analysis of properties of variational equations, like the geodesic
equation~\cite{Man03}, the minimal surfaces equation, etc..  Unfortunately, the
development of specific examples is still incomplete; we hope to be able to
produce them in a relatively short time.

\subsection*{Acknowledgements}

We would like to thank A.~M. Verbovetsky for many stimulating
discussions. We would also like to thank G.~De Cecco, L. Fatibene,
J. Jany\v ska, M. Palese.

We thank the city council of S. Stefano del Sole (AV, Italy), where this
research has been partially developed during the summer course `Diffiety
School: a school on the geometry of PDE', July 2002.

This research has been partially supported by EPSRC (England), GNSAGA of
Istituto Nazionale di Alta Matematica (Italy), Istituto Italiano per gli Studi
Filosofici, King's College (London), Ministero per l'Istruzione, l'Universit\`a
e la Ricerca (Italy), Universities of Lecce and Salerno.

\end{document}